\documentclass[12pt]{article}
\usepackage{amsfonts}
\usepackage{amssymb}
\usepackage{latexsym}
\title{The consistency strength of NFUB\\
Preliminary draft}
\author{Robert M. Solovay\thanks{PartiPart of this work was completed
during the 1997 Elsag-Bailey --- I. S. I. Foundation research meeting
on quantum computation.}\\
PO Box 5949\\
Eugene OR 97405\\
solovay@math.berkeley.edu}
\newcommand{\Beth}{\mbox{Beth}}
\newcommand{\x}{200}
\newcommand{\y}{200}
\newtheorem{theorem}{Theorem}[section]
\newtheorem{lemma}[theorem]{Lemma}
\newtheorem{proposition}[theorem]{Proposition}
\newtheorem{definition}{Definition}[section]
\begin{document}
\maketitle
\begin{abstract}
We show that the consistency strength of the system NFUB, recently
introduced by Randall Holmes, is precisely that of $ZFC- + \mbox{
``There is a weakly compact cardinal''.}$
\end{abstract}
\tableofcontents
\section{Introduction}

\subsection{NF and some variants.}

The unorthodox system of set-theory, New Foundations, [or NF for short]
was introduced by Quine. His approach for blocking the paradoxes of
naive set theory was to introduce a stratification condition in the
comprehension axiom. It is still an open problem whether NF can be
proved consistent relative to one of the usual flavors of set theory
such as ZFC. Indeed there is no  known proof of the consistency of NF
even relative to one of the standard large cardinal axioms.

That the system NF has some rather counterintuitive properties was
shown by Specker. He proved among other things that the axiom of choice
was refutable in NF and [as a corollary] that the axiom of infinity
was a theorem of NF.

Subsequently, Jensen introduced NFU, the [slight?] variant of NF in which
the axiom of extensionality was weakened to allow elements which are
not sets. Jensen was able to prove that this variant was consistent,
and that it was compatible with the axioms of choice and infinity. I
shall differ from Jensen's terminology by taking the axioms of choice
and infinity as two of the axioms of NFU. I shall also suppose that an
ordered pair opertion is introduced as one of the primitives of the
system, and that [for purposes of stratification] the type of the
ordered pair is the same as that of its members. [The usual Kuratowski
ordered pair has type two more than the types of its constituents and
is therefore less suitable in the context of NFU.]

The theory NFU has a universal class, $V$,  which contains all of its
subsets. One would suspect that a contradiction is near. But disaster
is escaped because the obvious map of $V$ onto the set of all
one-element sets can not be proved to exist because of failure of
stratification. We say that a set $S$ is {\em Cantorian}
if there is a bijection of $S$ with the set $USC(S)$ consisting of all
the singletons whose members lie in $S$. A set is {\em strongly
Cantorian} if the obvious map [which sends $x$ to $\{x\}$] provides a
bijection of $S$ with $USC(S)$. Thus, in the NFU context, the effect
of the Russell paradox is to show that $V$ is not Cantorian.

Holmes considered the system NFUA which is obtained from NFU by
adjoining the axiom that ``Every Cantorian set is strongly
Cantorian''. Correspondence with Holmes prompted me to work out the
precise consistency strength of NFUA [in work which is as yet
unpublished]. The theory NFUA is equiconsistent with the theory obtained
from ZFC by adjoining [for each positive integer $n$] an axiom which
asserts the existence of an $n$-Mahlo cardinal.

We caution the reader about a subtle point. Naively, this looks the
same as the theory which is obtained from ZFC by adjoining a single
axiom that asserts ``For every integer $n$, there is an $n$-Mahlo
cardinal.'' But in fact, this latter theory has strictly greater
consistency strength than the theory involving an infinite list of
axioms which we mentiond in the preceding paragraph.

\subsection{NFUB.}

Holmes also introduced a stronger theory NFUB, which we now describe.
In NFU we can develop a theory 
of ordinals in the spirit of Whitehead and Russell so that an ordinal
consists of the class of all well-orderings which are order-isomorphic
to a given well-ordering. An ordinal is Cantorian if the underlying
sets of the well-orderings which are its members are all Cantorian. It
is easy to see that the Cantorian ordinals [in the theory NFUA] form
an initial segment of the ordinals which is not represented by a set.
[That the collection of Cantorian ordinals does not form a set is a
variant of the Burali-Forti paradox.]

Let us say that a subcollection $S$ of the Cantorian ordinals is {\em
coded} if 
there is some set $s$ whose members among the Cantorian ordinals are
precisely the members of $S$. The system NFUB is obtained from the
system NFUA by adding an axiom schema which asserts that any
subcollection of the Cantorian ordinals which is definable by a
formula of the language of NFUB [possibly unstratified and possibly
with parameters] is coded by 
some set.

Holmes showed that the existence of a measurble cardinal implies the
consistency of NFUB. Prior to my work on the problem, it was open if a
Ramsey cardinal implied the consistency of NFUB and whether or not
NFUB implied the existence of $0^{\#}$. 

It follows from our main theorem that [in $ZFC$] a weakly compact
cardinal implies the consistency of NFUB and that NFUB does not prove
``$0^{\#}$ exists''.

Our main result is the following:
\begin{theorem}
\label{main_theorem}
The following theories are equiconsistent:
\begin{enumerate}
\item NFUB.
\item $ZFC-$ + ``There is a weakly compact cardinal'',
\end{enumerate}
\end{theorem}

Remarks:
\begin{enumerate}
\item $ZFC-$ is the theory consisting of all the axioms of $ZFC$
except the power set axiom. All the usual formulations of the notion
of ``weakly compact'' remain equivalent under $ZFC-$. The formulation
of ``weak compactnes'' that we shall actually use is: $\kappa$ is
weakly compact iff $\kappa$ is strongly
inaccessible and every $\kappa$-tree has a branch.

\item Officially, our metatheory is $ZFC-$\ \  [or equivalently,
$2^{nd}$-order number theory]. The reader who is familiar with turning
model-theoretic consistency proofs into syntactic ones will have no
trouble in formalizing our arguments in Peano arithmetic [or, indeed,
in primitive recursive arithmetic].

\end{enumerate}

The remainder of this paper is organized as follows. Section 2 will
derive the consistency of $ZFC-$ plus a weakly compact from that of
NFUB. This proof is a very slight extension of earlier work of Randall
Holmes. The proof of the other direction is much more difficult.
Section 3 sets up the problem by showing that it suffices to construct
a model of ZFC equipped with an automorphism and having various
desirable properties. Section 4 outlines the transfinite construction
of such a model, with the key difficulty of how to execute the
successor step postponed to Section 5. In section 5, we show how to
handle the successor step. Here we must finally exploit the
weakly compact cardinal.

\section{Getting weakly compact cardinals}

\subsection{Outline of the proof}

Our construction of a model of $ZFC-$ + ``There is a weakly compact
cardinal'' depends heavily on earlier work of Holmes which we shall
have to review. First, every model of NFU has associated to it a
certain model of $ZFC-$, $Z$, whose construction we will review. An
important ingredient of the structure of $Z$ is a certain endomorphism
$T$ which we shall also recall.

There is an appropriate notion of Cantorian for elements of $Z$: those
elements fixed by $T$. The key axiom  schema of NFUB implies that
the Cantorian elements of $Z$ are the sets of a certain model of KM,
which we dub the {\em canonical model}. [KM is the variant of
class-set theory exposed in the appendix to Kelley's {\em General
Topology}.]  As we shall explain below, it
makes sense to ask if the class of ordinals, $OR$, of a model of KM is
weakly compact. We show that for the canonical model, $OR$ is indeed
weakly compact.

Given any model of KM, there is associated another model of KM in
which, in an appropriate sense, $V=L$ holds. If $OR$ was weakly compact in
the original KM model, it will remain weakly compact in the new $L$-like
KM model.

Finally, it will be easy to turn this $L$-like KM model into a model of
$ZFC-$ + ``there is a weakly compact cardinal''.

Thus ends our outline. We turn to the details.

\subsection{Cardinals}

Until further notice, we are working in NFU.

As we indicated in the introduction, the usual treatment of cardinals
and ordinals in NFU is in the spirit of Russell and Whitehead.

\begin{definition}
Let $X, Y$ be sets. Then $X$ is {\em equipotent} with $Y$ iff there is
a bijection mapping $X$ onto $Y$.

Let $X$ be a set. Then 
\[ card(X) = \{Y\mid \mbox{ $Y$ is equipotent with $X$}\}\]
\end{definition}

$\lambda$ is a cardinal iff there is a set $X$ such that $\lambda =
card(X)$.

$CARD = \{\lambda\mid \lambda \mbox{ is a cardinal}\}$.

\subsection{The $T$ operation on cardinals.}

As we recalled in the introduction, if $X$ is a set, 
\[ USC(X) = \{ \{x\}\mid x \in X\}\]

This operation makes sense in NFU. [I.e., the definition of $USC(X)$
is {\em stratified}.]

\begin{definition}
Let $\kappa$ be a cardinal. Then $T(\kappa) = card(USC(X))$ for some
[any] $X \in \kappa$.
\end{definition}

Caution: The map on $CARD$ which sends $\kappa$ to $T(\kappa)$ is {\em not}
given by a set.

As usual, $V = \{x\mid x=x\}$. We set $\kappa_0 = card(V)$ and define
$\kappa_n$ [by induction on $n$ in the metatheory] by setting
$\kappa_{n+1} = T(\kappa_n)$.

We have:
\[ \kappa_0 > \kappa_1 > \kappa_2 > \ldots \]

\subsection{Ordinals}

The treatment of ordinals is quite analogous to that of cardinals.

\begin{definition}

Let $X$ be a set and $R$ a binary relation on $X$. {\em $R$} is a
linear ordering of $X$ if:
\begin{enumerate}

\item $R$ is transitive.

\item For any $x,y \in X$, exactly one of $xRy$, $yRx$, $x=y$ holds.
\end{enumerate}
\end{definition}

Comment: Thus our linear orderings are {\em strict}.

Let $X$ be a set and $R$ a linear ordering on $X$. Then $otp(\langle X,
R \rangle) = \{\langle X_1, R_1 \rangle\mid \mbox{$R_1$ is a linear ordering
of $X_1$ and $\langle X_1, R_1\rangle$ is order isomorphic to $\langle
X , R\rangle$}\}$.

\begin{definition}
Let $X$ be a set and $R$ a linear ordering of $X$. Then $R$ is a
well-ordering of $X$ iff for every non-empty $U \subseteq X$ there is a
$u \in U$ such that $(\forall v \in U) (u R v \mbox{ or } u=v)$.
\end{definition}

\begin{definition}
$\lambda$ is an {\em ordinal} if $\lambda = otp(\langle X, R \rangle)$
for some well-ordering $\langle X, R \rangle$.

$OR = \{\lambda \mid \lambda \mbox{ is an ordinal}\}$

\end{definition}

Since we have the axiom of choice we can identify each cardinal with
the corresponding initial ordinal. That is, if $\gamma = card(X)$,
then letting $R$ be a well-ordering of $X$ which is ``as short as
possible'', we identify $\gamma$ with th ordinal $otp(\langle X, R
\rangle)$. It is easy to check this doesn't depend on the choices of
$X$ and $R$.

\subsection{The $T$ operation on ordinals}

Let $\lambda$ be an ordinal. We define $T(\lambda)$ as follows: Let
$\lambda = otp(\langle X, R \rangle)$. Let $X_1 = USC(X)$. Define $R_1
= \{ \langle \{x\}, \{y\}\rangle\mid \langle x , y\rangle \in R\}$.

Then $R_1$ is a well-ordering and we set $T(\lambda) = otp(\langle X_1,
R_1 \rangle)$.

It is easy to check that the identification of cardinals with initial
ordinals is compatible with the $T$ operations in $CARD$ and $OR$.

Suppose that $\langle X, R \rangle$ is a well-ordering and that
$\lambda = otp(\langle X , R \rangle)$. We define a new well-ordering
as follows:

The underlying set $X^\star$ will consist the set of ordinals less
than $\lambda$. The ordering $R^\star$ will be the restriction to
$X^\star$ of the usual ordering on $OR$.

In orthodox set-theory, $\langle X^\star, R^\star \rangle$ would be
order isomorhic to $\langle X , R \rangle$. However, in NFU we have:

\begin{proposition}
The order type of $\langle X^\star, R^\star \rangle$ is
$T^2(\lambda)$.
[I. e., it's $T(T(\lambda))$.]
\end{proposition}

In this draft, we won't give the proof. A proof can certainly be found
in Holmes' forthcoming book on NFU.

A closely related fact is the following:
\begin{proposition}
\label{OR_length}
The order type of $OR$ [equipped with the usual ordering] is ${\kappa_2}^+$.
\end{proposition}

Here, we are identifying cardinals with initial ordinals and writing
$\lambda^+$ for the least cardinal greater than the cardinal
$\lambda$. [Of course, if $\lambda = \kappa_0$, then $\lambda^+$ is
undefined.] 

\subsection{The set $Z$.}
\label{z_dfn}

We continue to review Holmes' work. It turns out that there is a
natural model $Z$ of $ZFC-$ associated to any model of NFU. The
construction of $Z$ is just the adaptation to the NFU context of a
familiar construction 
used, for example, to get a model of $ZFC-$ from a model of
$2^{nd}$-order number theory.

We shall differ slightly in the details of our development of $Z$ from
the treatment in Holmes' text. But the two treatments are completely
equivalent and yield isomorphic versions of $Z$.

Let $X$ be a set and $R$ a binary relation on $X$. We list various
properties that the relation $R$ can have.

\begin{enumerate}
\item $R$ is {\em extensional} if whenever $x$ and $y$ are two
distinct members of $X$, then there is a $z \in X$ such that $z R x
\not\equiv z R y$.

\item $R$ is {\em well-founded} if whenever $U$ is a non-empty subset of
$X$, there is an element $u \in U$ which is $R$-least in $U$ in the
sense 
that for no $v \in U$ do we have $v R u$.

\item $R$ is {\em topped} if there is an element $t \in X$ [the {\em top}]
such that 
for any $x \in X$ there is a finite sequence $x_0, \ldots, x_n$ [$n$
can be $0$] with $x_0 = x$, $x_n = t$ and $x_i R x_{i+1}$ for $0 \leq
i < n$.
\end{enumerate}

If $X$ is a set and $R$ is a binary relation on $X$, then we let
$iso(X,R)$ be the set of all pairs $\langle X', R'\rangle$ such that
$R'$ is a binary relation on $X'$ and $\langle X', R'\rangle$ is
isomorphic to $\langle X, R \rangle$.

We are now in a position to define $Z$. It consists of all such $iso
(X,R)$ where $R$ is an extensional, well-founded, topped relation on
$X$.

We remark that the top of a well-founded topped relation is unique.

To explain the intuition behind the definition of $Z$ it helps to
leave NFU for the moment and work in $ZFC$. If $\langle X , R \rangle$
is an extensional well-founded relation, then the Mostowski collapse
theorem provides a transitive set $T$ and an isomorphism $\psi$ between
$\langle X , R\rangle$ and $\langle T, \epsilon_T \rangle$. [Here
$\epsilon_T$ is the restriction of the usual epsilon relation to $T$.]

Moreover, the set $T$ and the map $\psi$ are uniquely determined by
$X$ and $R$.

Now if $\langle X,R\rangle$ is extensional, well-founded and topped
[with top $t$], then $\langle X, R\rangle$ can be viewed as a code for
$\psi(t)$. 

\subsection{The binary relation $E$ on $Z$}
\label{e_dfn}

Still following Holmes, we now define an ``$\epsilon$-relation'' on $Z$,
which we dub $E$.

We define it in terms of representives. [It should, of course, be
checked that the definition does not depend on the choice of
representatives.]

So let $z_1$, $z_2$ be elements of $Z$ and let $\langle X_i, R_i
\rangle$ be a member of $z_i$. Let $t_i$ be the top of $\langle X_i,
R_i \rangle$.

Then $z_1 E z_2$ iff there is a map $\psi$ mapping $X_1$ injectively
into $X_2$ such that:
\begin{enumerate}
\item For $x,y \in X_1$, $x R_1 y$ iff $\psi(x) R_2 \psi(y)$.
\item Range($\psi$) is transitive in $X_2$. That is, if $x \in X_1$,
$y \in X_2$ and $y R_2 \psi(x)$, then there is a $z \in X_1$ with $y =
\psi(z)$.
\item $\psi(t_1) R_2 t_2$.
\end{enumerate}

\subsection{The model $\langle Z , E\rangle$.}

To understand the model $\langle Z,E\rangle$, it helps to return once
again to $ZFC$. 

Recall that, for $\kappa$ an infinite cardinal,  $H_\kappa$ is the
collection of sets whose transitive 
closures have cardinality less than $\kappa$.

The model $H_{\kappa^+}$ is a model of $ZFC-$.
It is never a model of $ZFC$ since it has a largest
cardinal [namely $\kappa$].

Let $\kappa$ be an infinite cardinal.
Here is a second-order characterization of when a relational structure
$\langle X,R\rangle$ is isomorphic to $\langle H_{\kappa^+}, \in \rangle$.
\begin{enumerate}
\item $R$ is extensional and well-founded.
\item For every $x \in X$, the set $\{y\mid y Rx\}$ has cardinality
$\leq \kappa$.
\item Let $U \subseteq X$ with $card(U) \leq \kappa$. Then there is a
$u \in X$ such that $U = \{y \in X \mid y R u\}$.
\end{enumerate}

Let us refer to this second order characterization [involving the
parameter $\kappa$] as $\Phi(\kappa)$. Returning to the NFU context,
we have the following characterization of $\langle Z, E \rangle$ up to
canonical isomorphism:

\begin{proposition}
$\langle Z, E\rangle$ satisfies the second-order sentence $\Phi(\kappa_2)$.
\end{proposition}

We do not give the proof. The $\kappa_2$ appears for the same reason
that the ${\kappa_2}^+$ appears in Proposition~\ref{OR_length}.

As an immediate corollary, $\langle Z,E\rangle$ is a model of $ZFC-$.

Hence it has its own notion of ordinal: the usual von Neumann
definition in which an ordinal $\gamma$ is equal to the set of ordinals
less than $\gamma$.

The ordinals we have previously used can be identified with elements
of $Z$ as follows. Let $\lambda$ be an ordinal and let $\langle X, R
\rangle$ be some representative of $\lambda$. We may assume that the
underlying set $X$ is not all of $V$; let $x$ be an element not in
$X$. We form a new binary structure as follows: $X' = X \cup \{x\}$.
$R'$ is the well-ordering of $X'$ that agrees with $R$ on $X$ and puts
$x$ ``at the end''. Then the structure $\langle X',R'\rangle$ is
extensional, well-founded and topped, and so determines an element of
$Z$. It is this element which we identify with $\lambda$. It is not
hard to see that this map of $OR$ into $Z$ is given by a set of NFU
and that it maps $OR$ onto the von Neumann ordinals of $Z$ in an
order-preserving fashion.

It is this writer's opinion that the best way to define set-theoretic
notions in NFU [at least in the presence of the axiom of choice] is to
reduce matters to working in $Z$ where stratification issues can for
the most part be avoided. When they can't be  avoided, they can be
reduced to a consideration of the $T$ operator to which we now turn.
\subsection{The $T$ operator on $Z$}

The definition of $T$ on elements of $Z$ is totally analogous to the
earlier definitions given for $OR$ and $CARD$.

Let $z \in Z$. Let $\langle X, R\rangle$ be a member of $z$. Let $X_1
= USC(X)$ and let $R_1 = \{\langle\{x\},\{y\}\rangle\mid x R y\}$
. Then
$\langle X_1, R_1 \rangle$ is well-founded, extensional and topped. We
set $T(z) = iso(X_1,R_1)$. It is clear that $T(z)$ does not depend on
the choice of representative $\langle X,R\rangle$.

We will need the following lemma in a moment:
\begin{lemma}
$2^{\kappa_1 }  $ is defined [and $ < \kappa_0$].
\end{lemma}

Proof sketch:

Define the ordinal $\eta$ by the requirement that $\Beth(\eta) \leq
\kappa_0$ and $\eta$ is largest so that this is true.

Applying $T$ we have $\Beth(T(\eta)) \leq \kappa_1$ and $T(\eta)$ is
largest such that this is true.

Now clearly (a) $T(\eta) \leq \eta$ [since $\kappa_1 < \kappa_0$] and
(b) $\eta$ is not Cantorian [since $\kappa_0$ is not Cantorian]. It
follows that $T(\eta) < \eta$.

$T(\eta)$
 has clearly the same residue class mod 3 as $\eta$ [since $T$
preserves all second order properties]. It follows that $T(\eta) + 3
\leq \eta$.

So we have $\kappa_1 \leq \Beth(T(\eta) + 1)$ so $2^{\kappa_1}
\leq \Beth(T(\eta) +2) < \Beth(\eta) \leq \kappa_0$. $\Box$

In order to understand the $T$ map on $Z$ we make the following
construction. Let $Z_1 = USC(Z)$. Let $E_1 = \{\langle
\{x\},\{y\}\rangle\mid xE y\}$. Then the pair $\langle Z_1,E_1\rangle$
satisfies the second order sentence $\Phi(\kappa_3)$ and so has
cardinality $2^{\kappa_3} \leq \kappa_2$. It is also clearly
well-founded and extensional. It follows that there is an isomorphism
$k$ mapping $\langle Z_1, E_1\rangle$  onto a transitive
set of the model $\langle Z, E \rangle$.

\begin{proposition}
Let $z \in Z$. Then $k(\{z\}) = T(z)$. 
\end{proposition}

It follows that the range of
the $T$ map on $Z$ is a set of the NFU model. The $T$ map itself is
definitely not given by a set of the NFU model.

The identification of the Russell-Whitehead ordinals with the
von Neumann ordinals of $Z$ is compatible with the $T$ operations on
the domain and range of the identification map.

\subsection{Cantorian elements of $Z$}

We say that an element of $Z$ is {\em Cantorian} if $T(z) = z$. The
main axiom of NFUA entails that the collection of Cantorian elements
of $Z$ is transitive. That is, if $z \in Z$ is Cantorian and $z_1 E z$
then $z_1$ is also Cantorian.

Using the axiom of choice, we can translate the main axiom scheme of
NFUB as
follows. Let $W$ be a subcollection of the elements of our NFUB model
that is definable by a formula of $\varphi$ of the language of NFUB,
possibly containing names for particular elements of our NFUB model.
[The formula $\varphi$ need not be stratified.] Then there is an
element $w$ of $Z$ which {\em codes} the intersection of $W$ with the
Cantorian elements of $Z$. That is, for any Cantorian element of $Z$,
say $c$, we have $c$ is in the collection $W$ iff $c E w$.

Holmes has shown that we get a model of $KM$ [Kelley-Morse set theory
including the global axiom of choice] as follows. The sets of the
model are the Cantorian elements of $Z$. The $\epsilon$-relation of
the model is the restriction of $E$. Finally, the classes of the model
are the subcollections of the Cantorian sets which are coded by some
element of $Z$ in the manner just described.

\subsection{Digression on weak compactness}

We first have to explicate  what it means to assert that the class of
ordinals of a model of $KM$ is weakly compact. But this is easy. The
various definitions of ``weakly compact'' are given by $\Pi^1_2$
formulas. We interpret the type $2$ variables as ranging over the
classes of the $KM$ model. We interpret the type $1$
variables as ranging over the sets of the $KM$ model.

There are many different characterizations of what it means for a
strongly inaccessible cardinal to be weakly compact. The notions that
will be important for us in this paper are:
\begin{enumerate}
\item $\kappa$ is $\Pi^1_1$-indescribable.

\item  $\kappa$ has the tree property. [There are actually two
versions of the tree property that we will need to consider. This will
be
spelled out in a moment.]
\end{enumerate}

We remark that the usual proofs that these formulations are equivalent
carry over to the $KM$ context without difficulty.

We turn to the precise formulation of the two versions of the tree
property we will need. In the following $\kappa$ is a strongly
inaccessible cardinal. At a first cut, the reader may take us as
working in ZFC though we will have occasion to apply these definitions
later in the context of $ZFC-$ and in the context of $KM$ [with
$\kappa$ replaced by the class of all ordinals].

A {\em tree} is a pair $\langle X, R \rangle$ where $X$ is a set, $R$
is a transitive binary relation on $X$ and whenever $a \in X$ and
$seg(a) = \{y \in X : y R a\}$ then $seg(a)$ is well-ordered by the
restriction of $R$.

If $\langle T, <_T \rangle$ is a tree and $t \in T$, then the rank of
$t$ [notation: $\rho(t)$] is the order type of the restriction of
$<_T$ to $seg(t)$.

$T$ is a $\kappa$-tree if:
\begin{enumerate}
\item For every $t \in T$, $\rho(t) < \kappa$.
\item If $\alpha < \kappa$, then let $T_\alpha = \{t \in T \mid
\rho(t) = \alpha\}$.

Then we require that for every $\alpha < \kappa$, $T_\alpha$ is
non-empty and has cardinality less than $\kappa$.
\end{enumerate}

A branch through a $\kappa$-tree $T$ is a subset $b$ of $T$ with the
following properties:
\begin{enumerate}
\item For every $\alpha < \kappa$, $b \cap T_\alpha$ has exactly one
member.
\item If $x \in b$ and $ y <_T x$, then $y \in b$.
\end{enumerate}

We can now give our first formulation of weak compactness: $\kappa$ is
weakly compact if $\kappa$ is strongly inaccessible and every
$\kappa$-tree has a branch.

It turns out that it is not necessary to consider all $\kappa$-trees.

A binary $\kappa$-tree is a set $S$ with the following properties:
\begin{enumerate}
\item The elements of $S$ are functions whose domain is an ordinal
less than $\kappa$ and whose range is included in $\{0,1\}$.
\item If $f \in S$ has domain $\alpha$ and $\beta < \alpha$, then the
restriction of $f$ to $\beta$ lies in $S$.
\item If $\alpha$ is an ordinal less than $\kappa$, then there is an
$f \in S$ with domain$(f) = \alpha$.
\end{enumerate}

Such an $S$ gives rise to a $\kappa$-tree iin our previous sense if we
take $<_S$ to be the restriction of $\subset$ to $S$.

What a branch amounts to for a binary $\kappa$-tree $S$ is a function $F:
\kappa \mapsto 2$ such that the restriction of $F$ to any ordinal less
than $\kappa$ lies in $S$.

Then an equivalent formulation of the notion of weak compactness is
that $\kappa$ is weakly compact iff $\kappa$ is strongly
inaccessible and every binary $\kappa$-tree has a branch.

\subsection{$OR$ is weakly compact}

It is now easy to show that in the canonical model of $KM$ [associated
to some given model of NFUB] $OR$ is weakly compact. Let then $S$ be a
class of the canonical model which gives a binary $OR$-tree. Let $s$
be an element of $Z$ which codes $S$. Then for arbitraily large
Cantorian ordinals $\alpha$, $\alpha$ is a Beth fixed point and $s
\cap V_\alpha$ is a binary $\alpha$-tree [in the sense of $Z$]. So
there must be a non-Cantorian $\alpha \in Z$ such that $Z$ thinks that
$\alpha$ is a Beth fixed point and that $s \cap V_
\alpha$ is a binary
$\alpha$-tree. [Otherwise, there would be a set of $Z$ consisting
precisely of the Cantorian ordinals, which is absurd.]

Let $b$ be an element of $s \cap V_\alpha$ of non-Cantorian rank. Then
it is easy to see that $b$ codes a branch through $S$ [consisting of
the restrictions of $b$ to Cantorian ordinals of $Z$]. Our proof that
$OR$ is weakly compact in the canonical model of $KM$ is complete.

\subsection{Getting a model of $ZFC-$}

The G\"{o}del $L$ construction applies to models of $KM$ as follows:
Given a well-ordering of $OR$, we can build a model of $V=L$ whose
ordinals have the order type of $R$. Say that a class is constructible
if it appears in some such model. Then if we take the sets of our new
model of KM to be the constructible sets of the old model, and the
classes of our new model to be the constructible classes of our old
model, we get a model of $KM + V=L$ in which every class is
constructible. 

The usual proof that if $\kappa$ is weakly compact, then it remains
weakly compact in $L$ works also in the present context. It gives that
$OR$ is still weakly compact in our $L$-like model of $KM$.

What we  have gained by this $L$-construction is a definable
well-ordering of the classes of our $KM$ model. In our original model,
we only had, a priori, a well-ordering of the class of all sets.

It is now routine [much as we built $Z$ in NFUB] to build a model of
$ZFC-\: +\; V=L$ + ``there is a weakly compact cardinal $\kappa$'' from
our $L$-like model of $KM$. This direction of the equiconsistency
proof is complete.

\section{Getting a model of NFUB}
\subsection{Initial preparations}
\label{m_dfn}

We are given a model of $ZFC-$ + ``there is a weakly compact
cardinal''. Say $\kappa$ is a weakly compact cardinal of the model. We
pass to the constructible sets of the original model. Then we still
have a model of $ZFC-$ and $\kappa$ remains weakly compact. Moreover
$V=L$ holds.

We may assume as well that $\kappa$ is the largest cardinal of the
model. For if not, let $\lambda = \kappa^+$ [in the sense of our
current model]. Then $L_\lambda$ is a model of $ZFC-$ + $V=L$ +
``There is a largest cardinal $\kappa$'' + ``$\kappa$ is weakly
compact''.

The upshot is that we may assume given a model $M$ of $ZFC-$ + $V=L$ +
``There is a largest cardinal $\kappa$ which is weakly compact''. We
must construct a model of NFUB. We shall reserve the symbols $M$ and
$\kappa$ for this model and this cardinal for the remainder of the paper.

By a class of $M$ we mean a subcollection of $M$ definable by a
formula of the language of set-theory [with the definition possibly
involving names for particular elements of $M$]. Our final model $N$
of NFUB will be a proper class model in the sense of $M$.

\subsection{Precise description of NFU}

We have been a bit lax in describing precisely the formal system NFU. We
now need to remedy this.

[Probably in the next draft of the paper, this material should be
placed much earlier.]

The language of NFU is a first-order language with no function symbols
and with the following five predicate symbols:
\begin{enumerate}
\item A binary predicate $=$. [This is the usual equality predicate.
The basic properties of $=$ are part of first order logic and need not
be specified in our axioms for NFU.]
\item A binary predicate $\in$. [$a \in b$ is read ``$a$ is a member
of $b$''.]
\item A unary predicate $S$. [$Sa$ is read ``$a$ is a set.'' Allowing
``urelements'' which are not sets is the key difference between NFU
and NF.]
\item A ternary precicate $P$. 

[The intuition  is as follows. Every pair of elements in
NFU determines a unique ordered pair.
 The relation $Pabc$ means ``$c$ is the ordered pair determined by 
$a$ and $b$''.]
\end{enumerate}

NFU will have three groups of axioms:

\subsubsection{Extensionality axioms}

\begin{enumerate}
\item If $a \in b$ then $Sb$.
\item Suppose that $Sa$, $Sb$ and $a \neq b$. Then there is a $c$ such
that $c \in a \not\equiv c \in b$.
\end{enumerate}

\subsubsection{Pairing axioms}
\begin{enumerate}

\item For every $a$ and $b$, there is exactly one $c$ such that $Pabc$.

\item If $Pabc$ and $Pa'b'c$, then $a = a'$ and $b = b'$.
\end{enumerate}

\subsubsection{Comprehension axioms}

The description of these is a bit more technical. Let $\varphi$ be a
formula of the language of NFUB. $\varphi$ is {\em stratified} if
there is a map $\sigma$ from the set of variables occurring in $\varphi$
[either free or bound] into $\omega$ [the set of non-negative
integers] such that:
\begin{enumerate}
\item If $v=w$ is a subformula of $\varphi$ then $\sigma(v) = \sigma(w)$.
\item If $Puvw$ is a subformula of $\varphi$ then $\sigma(u) = \sigma(v) =
\sigma(w)$.
\item If $v \in w$ is a subformula of $\varphi$ then $\sigma(w) =
\sigma(v) + 1$.
\end{enumerate}

Let now $\phi$ be a stratfied formula whose free variables are
included among $v_0, \ldots, v_n$.
Then the following is an axiom
of NFU:

$\forall v_1, \dots, v_n \exists v_{n+1} \forall v_0 [v_0 \in v_{n+1}
\iff \varphi]$.

\subsection{A procedure for getting models of NFU}
\label{q_dfn}

We describe a known procedure for getting models of NFU. The
conditions we impose on the starting model are far too stringent.
There is no need to require that $V=L$ holds in the model, or that the
model is a model of full $ZFC$. Moreover the requirement that $j$ be
an automorphism can be considerably relaxed. Nevertheless, we will
have these conditions in our applications and they serve to simplify
the discussion.

We consider the following situation:
\begin{enumerate}
\item $\langle N, \in_N\rangle$ is a model of $ZFC\ + V=L$.
\item $j$ is an automorphsim of the model $\langle N, \in_N\rangle$.
$j$ is not the identity.

It follows, [say since $V=L$ holds in $N$], that $j$ moves some
ordinal $\eta$. In fact, $j$ must move some Beth fixed point. [Indeed,
$j$ moves the $\eta^{th}$ Beth fixed point.]

\item $\gamma$ is a Beth fixed point of the model $N$ moved by $j$.
\end{enumerate}

Without loss of generality, we can assume that $j(\gamma) > \gamma$.
[Otherwise, simply replace $j$ by $j^{-1}$.]

{}From this data, we are going to define a model of NFU. We shall view
this procedure as ``well-known'' and not carry out the verification
that the model we describe is, indeed, a model of NFU.

We shall call the model we are constructing $Q$.

\begin{enumerate}
\item The underlying set of $Q$, $|Q|$ is just $L_\gamma$.

\item For $i \in \omega$, we define the ordinal $\gamma_i$ of the
model $N$ as follows:
\begin{eqnarray*}
\gamma_0 & = & \gamma;\\
\gamma_{n+1} & = & j^{-1}(\gamma_n).
\end{eqnarray*}

$Sx$ holds iff $x \subseteq L_{\gamma_1}$.
\item $=$ is the ordinary equality. Thus $x=y$ holds in $Q$ iff $x=y$.

\item However, the definition of $\in_Q$ is non-standard. $x \in_Q y$
 iff $Sy$  and $x \in_N j(y)$.

\item  $Pxyz$  holds in $Q$ iff $z$ is, in $L_\gamma$, the usual
Kuratowski ordered pair of $x$ and $y$. [That is $z = \{\{x\},\{x,y\}\}$.]
\end{enumerate}

This completes the specification of the model $Q$. As we have already
mentioned, $Q$ is a model of $NFU$. It turns out that the $T$
operation of $Q$ is essentially identical with $j^{-1}$. The relevant
facts will be recalled in the next subsection.

Remark: Let us write $\langle x,y\rangle$ for the unique $z$ such that
$Pxyz$ holds in $Q$. We caution the reader that, despite its origin
using the Kuratowski ordered pair in $L_\gamma$, $\langle x, y\rangle$
is {\em not} the Kuratowski ordered pair in the sense of the model
$Q$.

\subsection{$T$ vs. $j^{-1}$}

Let $Z$ be as defined in Section~\ref{z_dfn} [with repect to the model
$Q$].

Let $Z^\star = \{x \in |Q| \mid x \in_Q Z\}$. Here the definition takes
place in the model $N$ [using the parameter $j(Z)$].

Similarly, let $E^\star = \{\langle x, y\rangle \mid xEy\mbox{ holds
in }Q\}$.
Again, the definition can be given in $N$ using the parameter $j(E)$.

So, in $N$,  $E^\star$ is a binary relation on $Z^\star$. Intuitively,
this is $N$'s copy of the structure $\langle Z, E \rangle$ of $Q$.

In fact, $Z^\star$ can be given a direct description quite analogous
to the definition of $Z$ in $Q$. In place of considering binary relations on
subsets of $V$, we consider binary relations whose underlying set is a
subset of $L_{\gamma_2}$. Thus in our present context, for $X
\subseteq L_{\gamma_2}$ and $R$ a binary relation on $X$, we let
$iso(X,R)$ consist of all pairs $\langle X', R' \rangle$ where:
\begin{enumerate}
\item $X' \subseteq L_{\gamma_2}$;
\item $R'$ is a binary relation on $X'$;
\item The structure $\langle X', R' \rangle$ is isomorphic to the
structure $\langle X, R \rangle$.
\end{enumerate}
[Of course, this definition takes place within the model $N$.]

Finally $Z^\star$ consists of all sets of the form $iso(X,R)$ where:
\begin{enumerate}
\item $X \subseteq L_{\gamma_2}$;
\item $R$ is a binary relation on $X$ which is well-founded,
extensional, and topped.
\end{enumerate}
[Again, this definition takes place in $N$.]

The definition of $E^\star$ is the obvious analogue, in the present
context, of the prior definition of $E$. (Cf. section~\ref{e_dfn})

An immediate corollary of our description of $\langle Z^\star, E^\star
\rangle$ is that this structure is isomorphic to $L_{{\gamma_2}^+}$
[with its usual $\epsilon$-relation].

Let us use $k$ to describe this isomorphism. Thus if $z \in_Q Z$,
$k(z)$ is the element of $L_{{\gamma_2}^+}$ that corresponds under the
isomorphism just described.

The following proposition now explains the sense in which $T$ can be
identified with $j^{-1}$.

\begin{proposition}
Let $z \in_Q Z$. Then $k(T(z)) = j^{-1}(k(z)))$.
\end{proposition}

For the moment, I am taking this proposition as ``well-known''. I may
include a proof in a later draft of this paper.

\subsection{Criteria for $Q$ to be a model of NFUB}
\label{criteria}

First, here is a sufficient criterion for $Q$ to be a model of NFUA.

Criterion 1: Suppose that $N$, $j$, $\gamma$ are as in
section~\ref{q_dfn}. Suppose 
further that whenever $j(\alpha) = \alpha$ [for $\alpha \in OR$] and
$\beta < \alpha$ then $j(\beta) = \beta$. Then $Q$ is a model of NFUA.

This is immediate from the identification of $T$ with $j^{-1}$.

We call the Cantorian elements of $N$
those elements fixed by $j$. If $N$ satisfies criterion $1$, they will
form an initial segment of $N$. We let $C$ denote the collection of
Cantorian elements of $N$.

Let $S$ be a subcollection of the Cantorian elements of $N$. We say
that $S$ is coded by the element $s$ of $N$ if whenever $x \in C$,
then $(x \in S) \iff (x  \in_N s)$.

Let $N^\star$ be the structure $\langle N; \in_N, j\rangle$. We have
the obvious notion of a class of $N^\star$. This is a subset of $N$
definable by a formula of the language appropriate to $N^\star$
[possibly with parameters from $N$.] It is evident that any class of
$Q$ is a class of $N^\star$. Hence the following criterion is easy to
verify.

Criterion 2: Suppose that $N$ satisfies criterion 1, and that whenever
$W$ is a class of $N^\star$, then the intersection of $W$ with $C$ is
coded by some element $w$ of $N$. Then $Q$ is a model of $NFUB$.

\subsection{Some promises and their consequences}
\label{promises}

We are going to describe properties of a model $N_\infty$ and argue
that if we can construct a model with these properties then the
construction of section~\ref{q_dfn} will yield a model of NFUB.

Recall the model $M$ and its cardinal $\kappa$ that we introduced in
section~\ref{m_dfn}. 

The model $N_\infty$ will have the following properties:

\begin{enumerate}
\item It will be a class-sized model of $ZFC + V=L$. That is, both the
``underlying set'' and the $\epsilon$-relation of $N_\infty$ will be
classes of $M$.
\item There will be an elementary embedding $i_\infty: L_\kappa
\mapsto N_\infty$ which maps $L_\kappa$ onto an initial segment of
$N_\infty$.

[$i_\infty$ will be a set of $M$.]

\item 
\label{automorphism}
There will be an automorphism $j_\infty$ of the model
$N_\infty$ [again given by a class of $M$].

The only elements of $N_\infty$ left fixed by $j_\infty$ are those in 
the range of $i_\infty$.

\item 
\label{coding}
Let $A$ be a subset of $L_\kappa$ lying in $M$. Then there is an
element $a \in N_\infty$ which codes $A$ in the sense that for all $x
\in L_\kappa$ we have $x \in A \iff i_\infty(x) \in_{N_\infty} a$.
\end{enumerate}

It follows first from item~\ref{coding} that $N_\infty$ is a proper
class of $M$ and hence, since range $i_\infty$ is a set of $M$ that
$j$ is not the identity. Picking some Beth fixed point $\gamma$ of
$N_\infty$ which is moved by $j_\infty$, we are in position to apply
the construction of section~\ref{q_dfn} to get a model $Q$ of NFU.
Applying criterion 1 of section~\ref{criteria} and
item~\ref{automorphism} of the list of properties of $N_\infty$ we see
that $Q$ is a model of NFUA.

It remains to see that $Q$ is a model of NFUB. We seek to apply
criterion 2. First note that we can identify the Cantorian elements of
$N_\infty$ with $L_\kappa$ via the map $i_\infty$.

Because all the elements of $N_\infty^\star$ are classes of $M$ and
$M$ is a model of $ZFC-$, the intersection of any class of
$N_\infty^\star$ with the cantorian elements of $N_\infty$ will
correspond to a subset $S$ of $L_\kappa$ which appears in $M$. But
then item~\ref{coding} of the list of properties of $N_\infty$ will
guarantee that $S$ is coded by some element of $N_\infty$.

So to complete the construction of a model $Q$ of NFUB, it remains to
construct a model $N_\infty$ with the stated properties.

\section{The transfinite construction}

\subsection{Introduction to this section}

The model $N_\infty$ will be constructed by a transfinite
construction, carried out within $M$,  whose stages are indexed by
the ordinals of $M$. The various stages will have various additional
components [in addition to being models of $V=L$]. We introduce a
category of $A$-models, so each of the stages will be such an
$A$-model. 

We let $\langle S_\alpha \mid \alpha \in {OR}^M \rangle$ be a
listing of the subsets of $L_\kappa$, lying in $M$, in order of
construction. [I.e., as ordered by $<_L$.] The construction of
$N_{\alpha + 1}$ will be devoted to insuring that $S_\alpha$ is coded in
the final model. 

For limit ordinals $\lambda$, $N_\lambda$ will be
constructed by a direct limit proces; the final model $N_\infty$ will
also be constructed by a direct limit process. This is
straightforward, but we choose to carefully present the direct limit
process below. 

We remark that for $\alpha \in {OR}^M$, the model
$N_\alpha$ will be a set of $M$ [and in fact have cardinality
$\kappa$].

\subsection{The category of $A$-models: Objects}

The following definition takes place within the model $M$.

An $A$-model consists of the following data:
\begin{enumerate}
\item A model $N$ of $ZFC + V=L$ of cardinality $\kappa$;

\item An elementary embedding $i: L_\kappa \mapsto N$;

\item An automorphism $j$ of $N$.
\end{enumerate}

These data are subject to the following requirements:
\begin{enumerate}
\item $i$ maps $L_\kappa$ onto an initial segment of $N$  [with
respect to the ordering that is $N$'s version of $<_L$].
\item The  points of $N$ left fixed by $j$ are precisely those in the range
of $i$.
\end{enumerate}

Note that there is the trivial $A$-model in which $N = L_\kappa$ and
$i$ and $j$ are identity maps.

\subsection{The category of $A$-models: Maps}

Let ${\cal N} =\langle N, i, j \rangle$ and ${\cal{N}}' =\langle N',
i', j'\rangle$ be $A$-models. 
A map of $A$-models, $\pi: \cal{N} \mapsto {\cal{N}}'$ is an
elementary embedding $\pi: \langle N, \in_N\rangle \mapsto \langle N',
\in_{N'}\rangle$ such that: 
\begin{enumerate}
\item The map $\pi$ should respect the embeddings $i$ and $i'$. That
is, the diagram indicated in figure~1 should commute.
\begin{figure}
\begin{center}
\begin{picture}(\x,\y)
\put(20,180){$L_\kappa$}
\put(25,160){\vector(0,-1){120}}
\put(40,160){\vector(1,-1){125}}
\put(20,20){$N$}
\put(168,20){$N'$}
\put(35,23){\vector(1,0){126}}
\put(15,100){$i$}
\put(110,100){$i'$}
\put(94,28){$\pi$}
\end{picture}
\end{center}
\caption{}
\end{figure}
\item The map $\pi$ should appropriately intertwine the automorphisms
$j$ and $j'$. That is, the diagram indicated in figure~2 should
commute.
\begin{figure}
\begin{center}
\begin{picture}(200,200)
\put(20,180){$N$}
\put(180,180){$N'$}
\put(20,20){$N$}
\put(180,20){$N'$}
\put(25,175){\vector(0,-1){142}}
\put (186,175){\vector(0,-1){142}}
\put(37,24){\vector(1,0){138}}
\put(37,184){\vector(1,0){138}}
\put(100,28){$\pi$}
\put(100,188){$\pi$}
\put(30,100){$j$}
\put(190,100){$j'$}
\end{picture}
\end{center}
\caption{}
\end{figure}
\end{enumerate}

It is easy to see that the identity map of an $A$-model is an
$A$-model map and that $A$-model maps are closed under composition.
I.e., the $A$-models form a category.

\subsection{Inductive requirements}

We continue to work within the model $M$.

Our construction will proceed in stages indexed by the ordinals [of
$M$]. Let $\lambda$ be such an ordinal. Here are the inductive
requirements that we will maintain before stage $\lambda$:
\begin{enumerate}
\item For each $\alpha < \lambda$ we will have defined an $A$-model,
${\cal N}_\alpha$.
\item If $\alpha_1 \leq \alpha_2 < \lambda$, we will have defined an
$A$-model map $\pi_{\alpha_1,\alpha_2}: {\cal N}_{\alpha_1} \mapsto
{\cal N}_{\alpha_2}$.
\item If $\alpha < \lambda$, then $\pi_{\alpha,\alpha}$ is the
identity map of ${\cal N}_\alpha$.
\item If $\alpha_1 \leq \alpha_2 \leq \alpha_3 < \lambda$, then the
map $\pi_{\alpha_1,\alpha_3}$ is equal to the composition
$\pi_{\alpha_2,\alpha_3}\pi_{\alpha_1,\alpha_2}$.
\item Suppose that $\alpha + 1 < \lambda$. Let $S_\alpha$ be the
$\alpha^{th}$ subset of $L_\kappa$ in order of construction. Then
there will be an element $s \in N_{\alpha + 1}$ that codes $S$ in the
sense described in section~\ref{promises}. That is, for any $x \in
L_\kappa$ we have $x \in S_\alpha \iff i_{\alpha+1}(x) \in_{N_{\alpha
+ 1}} s$.
\end{enumerate}

\subsection{Continuing the construction}

Suppose that we are at stage $\lambda$ and that our inductive
requirements hold at that stage. Here is what we do:
\begin{enumerate}
\item $\lambda = 0$. We take ${\cal N}_0$ to be the trivial $A$-model.
[So $N_0 = L_\kappa$ and $i_0$ and $j_0$ are identity maps.]
\item $\lambda = \alpha + 1$. We use the following lemma which will be
proved in section~\ref{atomic}:

\begin{lemma}
\label{main_lemma}
Let $S$ be a subset of $L_\kappa$ and let ${\cal N}$ be an $A$-model.
Then there is an $A$-model ${\cal N}'$ and an $A$-model map $\pi:
{\cal N} \mapsto {\cal N}'$ such that:
\begin{enumerate}
\item The map $\pi$ does not map $N$ onto $N'$.
\item There is an element $s \in N'$ which codes the set $S$.
\end{enumerate}
\end{lemma}

We apply this lemma in the obvious way with ${\cal N}_\alpha$ in the
role of ${\cal N}$ and $S_\alpha$ in the role of $S$.

We set ${\cal N}_{\alpha_ + 1}$ equal to the $A$-model ${\cal N}'$
provided by the lemma. and set $\pi_{\alpha, \alpha + 1}$ equal to the
map $\pi$ provided by the lemma. For other $\xi \leq \alpha + 1$ we
define $\pi_{\xi, \alpha + 1}$ in the unique way that maintains our
inductive requirements.

\item $\lambda$ is a limit ordinal. Then we take ${\cal N}_\lambda$ to
be the direct limit of the system of $A$-models  and $A$-model maps
already defined. This direct limit construction will be reviewed in
the next subsection. It provides all the maps needed to maintain our
inductive requirements.
\end{enumerate}

\subsection{The direct limit construction}

Suppose we are at stage $\lambda$ for  some limit ordinal $\lambda$. We
describe the construction of ${\cal N}_\lambda$ and the associated
maps.

Let $\alpha < \lambda$ and let $x \in N_\alpha$. We say that $x$ is
{\em original} if for no $\beta < \alpha$ and $y \in N_\beta$ do we
have: $\pi_{\beta,\alpha}(y) = x$. It is clear that for every $x \in
N_\alpha$ there is a $\beta \leq \alpha$ and an original $y \in
N_\beta$ such that $\pi_{\beta,\alpha}(y) = x$.

The underlying set of ${\cal N}_\lambda$ will consist of all pairs $\langle
\alpha, y\rangle$ where $\alpha < \lambda$ and $y \in N_\alpha$ is
original.

Of course, we set $\pi_{\lambda,\lambda}$ equal to the identity map. 

Let now $\alpha < \lambda$ and let $x \in N_\alpha$. We have to define
$\pi_{\alpha,\lambda}(x)$. Let $\beta \leq \alpha$ and $y \in N_\beta$
with $\pi_{\beta,\alpha}(y) = x$ and $y$ an original element of
$N_\beta$. [These requirements clearly uniquely determine $y$ and
$\beta$.] Then we set $\pi_{\alpha,\lambda}(x) = \langle \beta, y
\rangle$.

Let $\langle \alpha, y \rangle$ be an element of $N_\lambda$. We set
$j_\lambda(\langle \alpha, y \rangle) = \langle \alpha,
j_{\alpha}(y)\rangle$. [It must of course be checked that
$j_{\alpha}(y)$ is an original element of $N_\alpha$. This is not hard
to do.]

Let $R$ be one of the relations $=$ and $\in$. Let $x$ and $y$ be
elements of $N_\lambda$. We must determine whether or not $x R y$
holds in ${\cal N}_\lambda$. Here is the procedure. Find $\gamma <
\lambda$ such that there are $x'$ and $y'$ in $N_\gamma$ with
$\pi_{\gamma,\lambda}(x') = x$ and $\pi_{\gamma,\lambda}(y') = y$.
Then $xRy$ holds in ${\cal N}_\lambda$ iff $x' R y'$ holds in ${\cal
N}_\gamma$. [It must, of course, be checked that this procedure does
not depend on the choice of $\gamma$.]

Finally, we take $i_\lambda$  to be $\pi_{0,\lambda}$. 

The check that, with these definitions, ${\cal N}_\lambda$ is an
$A$-model and our inductive conditions continue to hold 
left to the reader.

We remark only that the reason that $N_\lambda$ has cardinality at
most $\kappa$ is that this is true of $\lambda$.

\subsection{Defining ${\cal N}_\infty$}

The definition of ${\cal N}_\infty$ is totally analogous to the direct
limit procedure of the preceding section. The only difference is that
we take the direct limit of the full system $\langle {\cal
N}_\alpha\mid \alpha \in {OR}^M\rangle$ so $N_\infty$ is a proper
class.

We discuss only the coding property. Let then $S \subseteq L_\kappa$
be a set in $M$. Then $S = S_\alpha$ for some ordinal $\alpha$ and
there is an element $s \in N_{\alpha + 1}$ that codes $S$. But then
$\pi_{\alpha + 1, \infty}(s)$ codes $S$ in $N_\infty$.

\section{Using the weakly compact cardinal}
\label{atomic}
\subsection{Introduction to this section}
It remains to prove Lemma~\ref{main_lemma}. Our proof of this is
rather mysterious in that for most of the proof we are engaged in
constructions having nothing to do with the main lemma, and at the
last minute we return to it and prove it. Perhaps the following
comments will dispel some of the mystery.

First, it turns out that the heart of the problem is the special case
when the inital $A$-model is the trivial one. In this case, the new
model ${\cal N}'$ will be generated [in an appropriate sense] from
$L_\kappa$ and a sequence of ``indiscernibles'' $\langle \xi_i\mid i
\in \mbox{\bf Z}\rangle$. [Here {\bf Z} is the set of {\em all}
integers, positive, negative, or zero.]

Everything boils down to determining a suitable EM-blueprint for the
indiscernibles. It's natural to use the partition properties associated
with a weakly compact, but since we have $\kappa$ many partitions to
tame, this doesn't seem to work. What does work is to imitate the
proof of the partition properties, which we briefly recall.

Let $[\kappa]^n$ be the set of $n$-element subsets of $\kappa$.
[Recall that $\kappa$ is equal to the set of ordinals less than
$\kappa$.] We are given a map $F:[\kappa]^n \mapsto 2$ and we seek to
find a homogeneous set for $F$ of size $\kappa$. The approach is to
build a tree $T$ so that on any branch $b$ through the tree $F$ does
not depend on its last coordinate. This allows us to reduce the
problem from $n$ to $n-1$. 

In our case, this approach of simplifying situations on the branches
of trees will still apply. We will, in fact, start with a length
$\omega$ process of building trees and branches. Eventually, this will
help us in building EM blueprints and thereby models with
automorphisms.

\subsection{The basic module}
\label{basic-module}
The following lemma is where we make use of the fact that $\kappa$ is
weakly compact. It is indeed not difficult to deduce that $\kappa$ is
weakly compact if the lemma holds.

\begin{lemma}
\label{basic_module}
Suppose given the following data:
\begin{enumerate}
\item A subset of $\kappa$, $B$, of cardinality $\kappa$.
\item A $\kappa$-sequence of ordinals $\langle \gamma_i\mid i <
\kappa\rangle$ such that $0 < \gamma_i < \kappa$ for all $i < \kappa$.
\item A $\kappa$-sequence of functions $\langle F_i \mid i <
\kappa\rangle$ such that $F_i: \kappa \mapsto \gamma_i$.
\end{enumerate}

Then there is a subset $B'$ of $B$ of cardinality $\kappa$ such that
each $F_i$ is constant on a tail of $B'$. That is, for every $i <
\kappa$ there is an $\eta < \kappa$ and a $\xi < \gamma_i$ such that
$F_i(\alpha) = \xi$ whenever $\alpha \in B'$ and $\alpha \geq \eta$.
\end{lemma}

We have stated the lemma in a form most useful for applications. The
version when $B$ is taken to be $\kappa$ is obviously equivalent.

Proof: We define a $\kappa$-tree $T_1$ as follows. The nodes of $T_1$
consist of functions $h$ whose domain is an ordinal $\alpha < \kappa$
and such that $h(\beta) < \gamma_\beta$ for all $\beta < \alpha$.

Let $\langle b_\xi\mid \xi < \kappa\rangle$ be an enumeration of $B$
in increasing order. We will, by induction on $\xi$, assign a node of
$T_1$ to $b_\xi$. We will do this so that:
\begin{enumerate}
\item At most one ordinal is assigned to each node of $T_1$.
\item Suppose that $x$ and $y$ are nodes of $T_1$ with $x <_{T_1} y$.
Then no ordinal can be assigned to $y$ while the node $x$ is still unoccupied.
\end{enumerate}

So suppose that for all $\alpha < \xi$, $b_\alpha$ has been assigned
to some node of $T_1$. We show how to assign $b_\xi$. Define a
function $G_\xi:\kappa \mapsto \kappa$ as follows. $G_\xi(\alpha) =
F_\alpha(b_\xi)$. It is clear that for all $\alpha < \kappa$, the
restriction of $G_\xi$ to $\alpha$ is a node of $T_1$. We take
$\alpha$ as small as possible so that the restriction of $G_\xi$ to
$\alpha$ is currently unoccupied, and assign $b_\xi$ to this
restriction.

Let $T$ be the set of nodes of $T_1$ that are assigned some element of
$B$ in the procedure just described. It is clear that $B$ has
cardinality $\kappa$ and hence that it is a $\kappa$-tree. Let
$G:\kappa \mapsto \kappa$ give a branch through $T$. [It is here that we
are using the weak compactness of $\kappa$.] Let $B'$ be the set of
elements of $B$ which have been attached to nodes on this branch. Then
if $G(\alpha) = \xi$, then it is clear that $F_\alpha(\theta) = \xi$
for a tail of $\theta$'s in $B'$.

\subsection{A technical lemma}

Recall that throughout we are working in the model $M$ which is a
model of $ZFC$ + $V=L$ + ``There is a largest cardinal $\kappa$ which
is weakly compact''.

\begin{lemma}
Let $A \subseteq \kappa$. Let $\kappa < \alpha$. Then there is a
$\beta > \alpha$ such that $A \in L_\beta$ and $L_\beta$ is a model of
$ZFC-$ + $V=L$ + ``$\kappa$ is the largest cardinal''.
\end{lemma}

Proof: This is an easy consequence of the fact that $\kappa$ is
$\Pi^1_1$-indescribale. Note that we do not claim [and cannot claim]
that $\kappa$ is weakly compact in $L_\beta$.

\subsection{The length $\omega$ construction}
\label{length-omega}
As we mentioned in the introduction to this section, the construction
that follows is inspired by the constructions of finite length used to
prove that weakly compact cardinals have partition properties.

Let $A_0 \subseteq L_\kappa$. We are going to define the following [for
$i \in \omega$].
\begin{enumerate}
\item A subset $A_i \subseteq L_\kappa$.
\item A subset $B_i \subseteq \kappa$.

$B_0$ will be $\kappa$. $B_{i+1}$ will be a $\kappa$-sized subset of
$B_i$.

\item An ordinal $\xi_i > \kappa$
\end{enumerate}

Here is a description of the inductive construction of these objects.
First we handle the case $n=0$:
\begin{enumerate}
\item $A_0$ is given as the input to this construction.
\item $B_0 = \kappa$.
\item $\xi_0$ is the least ordinal greater than $\kappa$ such that
$A_0 \in L_{\xi_0}$ and $L_{\xi_0}$ is a model of $ZFC-$. [It follows
that $\xi_0$ has cardinality $\kappa$ and that $L_{\xi_0}$ thinks that
$\kappa$ is the largest cardinal.]
\end{enumerate}

Next we handle the case of $i=n+1$ assuming these objects have been
defined for $i =  n$:
\begin{enumerate}
\item $A_{n+1}$ is some simple encoding of the pair of sets $A_n$ and
$B_n$. To be definite, $A_{n+1} = [\{0\} \times A_n] \cup [\{1\}
\times B_n]$.
\item $B_{n+1}$ will be obtained from $B_n$ using
Lemma~\ref{basic_module}.

We fix the $L$-least enumeration of length $\kappa$ of the set of all
pairs $\langle F, \gamma\rangle$ such that $F \in L_{\xi_n}$, $F:
\kappa \mapsto \gamma$ and $0 < \gamma < \kappa$. We apply the lemma
to this enumeration and the set $B_n$, getting back a set which we
take for $B_{n+1}$. [Actually, the lemma only asserts that there {\em
exists} a set $B'$ with certain properties. We take $B_{n+1}$ to be
the $L$-least set satisfying the conclusions of the lemma.]

\item $\xi_{n+1}$ is the least ordinal $> \xi_n$ such that $A_{n+1} \in
L_{\xi_{n+1}}$ and $L_{\xi_{n+1}}$ is a model of $ZFC-$. It follows
that $\xi_{n+1}$ has cardinality $\kappa$ and that $L_{\xi_{n+1}}$
thinks that $\kappa$ is the largest cardinal.
\end{enumerate}

We set $\xi_\infty = \sup\{\xi_n\mid n \in \omega\}$.

The intuition is that the sets in $L_{\xi_\infty}$ are those over
which the length $\omega$ construction has control. Cf. the partition
result of the next subsection.

For the remainder of the proof, we keep the notation of this
subsection. We will presently exploit our freedom to choose the
starting set $A_0$.

\subsection{Partition theorems}
\label{partition}

Let $G : \kappa \mapsto \kappa$.

An increasing $n$-tuple of ordinals less than $\kappa$:
\[ \alpha_1 < \alpha_2 < \ldots < \alpha_n \]
is {\em $G$ spread apart} iff:
\begin{enumerate}
\item $G(0) < \alpha_1$;

\item For $1 < i \leq n$, $G(\alpha_{i-1}) < \alpha_i$.
\end{enumerate}

Recall that $[\kappa]^n$ is the set of size $n$ subsets from $\kappa$
[or what is much the same thing, the set of increasing $n$-tuples from
$\kappa$]. 

\begin{lemma}
\label{partition_lemma}
Let $\gamma$ be an ordinal with $0 < \gamma < \kappa$. Let $F:
[\kappa]^n \mapsto \gamma$ with $F \in L_{\xi_\infty}$. Then there is
an $m \in \omega$, a $G: \kappa \mapsto \kappa$ with $G \in
L_{\xi_\infty}$ and an $\eta < \gamma$ such that whenever $\alpha_1 <
\ldots < \alpha_n$ is an increasing $n$-tuple of ordinals from $B_m$
which is $G$ spread apart then
\[ F(\alpha_1, \ldots, \alpha_n) = \eta\]

Moreover, the value of $\eta$ does not depend on the choices of $m$
and $G$.
\end{lemma}

Proof: We first prove the claim of the first paragraph of the lemma
and proceed by induction on $n$. 

Consider first the case when $n = 1$. Pick $r$ large enough that $F
\in L_{\xi_r}$. Then the pair $\langle F, \gamma\rangle$ was considered when
constructing $B_{r+1}$> It follows that there is an $\eta < \gamma$
and an $\alpha < \kappa$ such that $F(\beta) = \eta$ whenever $\beta >
\alpha$ and $\beta \in B_{r+1}$. So it suffices to take $G$ the
function constantly equal to $\alpha$ and $m = r+1$.

We proceed to the inductive step where $n = m+1$ and we know the first
part of the lemma when $n = m$. [Here $n \geq 2$.] Let again $F \in
L_{\xi_r}$. Then clearly for every choice of $\alpha_1 < \ldots <
\alpha_m$, there is an ordinal $\beta$ and an ordinal $\eta < \gamma$
such that $F(\alpha_1, \ldots, \alpha_m, \beta^\star) = \eta$ whenever
$\beta^\star \in B_{r+1}$ and $\beta^\star > \beta$.

We can find a function $G^\star:\kappa \mapsto \kappa$ and a function
$H: [\kappa]^m \mapsto \gamma$ (both in $L_{\xi_{r+1}}$) which express the
dependence of $\beta$ and $\eta$ in the preceding paragraph on
$\alpha_1, \ldots, \alpha_m$. Namely, whenever 
\begin{enumerate}
\item $\alpha_1 < \ldots < \alpha_m < \kappa$;
\item $\kappa > \beta^\star > G^\star(\alpha_m)$;
\item $\beta^\star \in B_{r+1}$
\end{enumerate}

then $F(\alpha_1, \ldots, \alpha_m, \beta^\star) = H(\alpha_1, \ldots,
\alpha_m)$. 

We can now apply our inductive hypothesis to $H$ getting an integer $s
\in \omega$, a function $G_1: \kappa \mapsto \kappa$ in
$L_{\xi_\infty}$ and an $\eta < \gamma$ such that whenever
$\alpha_1 < \ldots \alpha_m$ lie in $B_s$ and are $G_1$ spread apart,
then $H(\alpha_1, \ldots, \alpha_m) = \eta$.

Define $G$ by $G(\alpha) = \max(G^\star(\alpha), G_1(\alpha))$. Let $s^\star
= \max(r+1,s)$. Then clearly whenever $\alpha_1, \ldots, \alpha_{m+1}$
lie in $B_{s^\star}$ and are $G$ spread apart, then \[F(\alpha_1,
\ldots, \alpha_{m+1}) = \eta.\] We have successfully completed the
inductive step.

We turn to the last paragraph of the lemma. Suppose that $m, G$ and
$\eta$ are such that whenever $\alpha_1, \ldots, \alpha_n$ lie in
$B_m$ and are $G$ spread apart, then $F(\alpha_1, \ldots, \alpha_n) =
\eta$.

Suppose further that $m', G'$ and
$\eta'$ are such that whenever $\alpha_1, \ldots, \alpha_n$ lie in
$B_{m'}$ and are $G'$ spread apart, then $F(\alpha_1, \ldots, \alpha_n) =
\eta'$. We must show that $\eta = \eta'$. 

This is not difficult. Let $m^\star = \max(m,m')$. Define $G^\star:
\kappa \mapsto \kappa$ by ${G^\star}(\alpha) = \max(G(\alpha),
G'(\alpha))$.

We can clearly find $\alpha_1, \ldots, \alpha_n$ in $B_{m^\star}$
which are $G$ spread apart. But then:
\[ \eta = F(\alpha_1, \ldots, \alpha_n) = \eta'. \]

\subsection{The ``ultrapower'' construction}
\label{ultrapower}
The word ``ultrapower'' is in quotes since what we do, though inspired
by the ultrapower construction, is somewhat different.

Throughout this section ${\cal A} = \langle A; R_1, \ldots R_k,
f_1,\ldots f_p\rangle$ is a first-order structure which is a member of
$L_{\xi_\infty}$. Thus $A$ is a set [necessarily of cardinality $\leq
\kappa$], $R_i$ is a relation on $A$ of arity $n_i$ and $f_i$ is an
operation of arity $m_i$. We let ${\cal L}$ be the first order
language appropriate to the similarity type of ${\cal A}$.

Our construction will give us a new model ${\cal A}^\star$ of the same
similarity type as ${\cal A}$ together with an elementary embedding
$d: {\cal A} \mapsto {\cal A}^\star$.

What is novel about our construction and differs from the usual
ultrapower construction is that there will also be a canonical
automorphism $k$ of ${\cal A}^\star$. The points of $A^\star$ which
are fixed by $k$ are precisely those in the range of the diagonal map
$d$.

\subsubsection{The ``measure space''}

We let $X$ be the set of all functions mapping {\bf Z} into $\kappa$.
[Any such function lies in $L_\kappa$, so the set $X$ is known to
$L_{\xi_\infty}$.] Here {\bf Z} is, of course, the set of integers
[positive, negative, or zero].

Let $f \in X$ and let $s$ be a finite subset of {\bf Z}. Then we write
$f[s]$ to indicate the restriction of $f$ to the set $s$.

\subsubsection{The class of functions}

In the ordinary definition of an ultrapower, we would consider all
functions from the measure space $X$ to the underlying set $A$ of our
target model. Here, however, we must put several restrictions on our
functions. 

We define a class ${\cal F}$ of functions mapping $X$ to $A$ as
follows. $f \in F$ iff:
\begin{enumerate}
\item $f$ is in the model $L_{\xi_\infty}$.
\item There is a finite subset $s$ of {\bf Z} such that the value of
$f(x)$ depends only on $x[s]$. [We say that $s$ is a support for the
map $f$.]
\end{enumerate}

Thus an $f \in {\cal F}$ can be described as the composition of three
maps:
\begin{enumerate}
\item The map which sends an element $x \in X$ to $x[s]$;

\item Let $s$ have $n$ elements: $s = \{s_1,\ldots, s_n\}$ where $s_1
< s_2 < \ldots s_n$. Then we have the map which sends x[s] to an
element of $\kappa^n$:$
 x[s] \mapsto \langle x(s_1), \ldots
x(s_n)\rangle$.
\item The final element of the composition is some function $F:
\kappa^n \mapsto A$ with $F \in L_{\xi_\infty}$.

N. B. The set $\kappa^n$ is the set of [not necessarily increasing]
$n$-tuples from $\kappa$. It should not be confused with the set
$[\kappa]^n$ of strictly increasing $n$-tuples from $\kappa$.
\end{enumerate}

\subsubsection{The ``ultrafilter''}

For each triple $\langle s, m,G\rangle$ where:
\begin{enumerate}
\item $s$ is a finite subset of {\bf Z};
\item $m \in \omega$;
\item $G: \kappa \mapsto \kappa$ lies in $L_{\xi_\infty}$
\end{enumerate}
we associate a subset $ A_{s,m,G} \subseteq X$ as follows:

An element $x \in X$ lies in $ A_{s,m,G}$ iff:
\begin{enumerate}
\item The function $x[s]$ is strictly increasing on its domain.
\item Let $s$ have $n$ elements. Let these elements, listed in
increasing order, be $s_1, \ldots, s_n$. Then $x(s_i) \in B_m$ for $1
\leq i \leq n$.
\item $x(s_1), \ldots, x(s_n)$ are $G$-spread apart.
\end{enumerate}

The reader should verify that the set $ A_{s,m,G}$ is non-empty.

These sets form the base of a filter in the following sense.
Let $\langle s_1,m_1,G_1\rangle$ 
and $\langle s_2,m_2,G_2\rangle$ be two triples of the sort just
discussed.

Set $s = s_1 \cup s_2$. Set $m = \max(m_1,m_2)$. Define $G:\kappa
\mapsto \kappa$ by setting $G(\alpha) =
\max(G_1(\alpha),G_2(\alpha))$.
Then:
\[ A_{s,m,G} \subseteq A_{s_1,m_1,G_1}\cap A_{s_2,m_2,G_2}.\]

The following terminology will make the analogy with the usual
ultrapower construction more transparent. Say that a subset of $X$ is
{\em measurable} if its characteristic function lies in the obvious
analogue of ${\cal F}$. [Make the same definition but replace reference
to $A$ by reference to $\{0,1\}$.] We
say that a measurable set has measure $1$ if it contains a set of the
form $A_{s,m,G}$. We say that a measurable set has measure $0$ if its
complement has measure $1$. It follows from the fact that sets of the
form $A_{s,m,G}$ form the base for a filter, that no set has
simultaneously measure $1$ and measure $0$. It is an immediate
consequence of the main result of Section~\ref{partition} that given a
measurable set $A \subseteq X$ precisely one of $A$ and $X-A$ has
measure $1$. In this way we have defined an ultrafilter on the Boolean
algebra of measurable sets. It is definitely not countably complete.

\subsubsection{Construction of the ``ultrapower'' ${\cal A}^\star$}

We omit many details since this now parallels the usual ultraproduct
construction. 

We first put an equivalence relation on the functions of ${\cal
F}$. Two such functions, say $f_1$ and $f_2$ are equivalent if
$\{x\mid f_1(x) = f_2(x)\}$ has measure $1$. $A^\star$ will consist of
all the equivalence classes of functions in ${\cal F}$ [for the
equialence relation just described].

The basic relations of ${\cal A}^\star$ are defined in terms of
representatives; it must be verified that the result is independent of
the choices made. Let $R_i$ be a $n$-ary relation of ${\cal A}$. The
analogue for ${\cal A}^\star$ ($R_i^\star$) is defined as follows:

$R_i^\star([f_1], \ldots, [f_n])$ holds in ${\cal A}^\star$ iff
$\{x\mid R_i(f_1(x), \ldots, f_n(x))\}$ has measure $1$.

Similarly, let $g$ be an $n$-ary operation of ${\cal A}$. The
corresponding operation of ${\cal A}^\star$, call it $g^\star$, is
defined thus:
$g^\star([f_1], \ldots, [f_n])$ is represented by the function:
\[ x \mapsto g(f_1(x), \ldots, f_n(x))\]

This completes our definition of the structure ${\cal A}^\star$. Note
that the resulting structure is clearly a {\em set} of $M$ [since the
functions we employ are taken from the {\em set} $L_{\xi_\infty}$; in
particular, the cardinality of the underlying set $A^\star$ is clearly
at most $\kappa$.

\subsubsection{The {\L}os theorem}
\begin{proposition}
Let $\varphi(v_1,\ldots v_n)$ be a formula of ${\cal L}$ [the language
appropriate to ${\cal A}$] having at most the indicated free
variables. Let $[f_1], \ldots, [f_n]$ be elements of $A^\star$.

Then $\varphi([f_1], \ldots, [f_n])$ holds in ${\cal A}^\star$ iff
\[\{x \in X \mid \varphi(f_1(x), \ldots, f_n(x)) \mbox{ holds in }
{\cal A}\}\] has measure $1$.
\end{proposition}

Proof: The usual proof applies without essential change.

\subsubsection{The diagonal map}

We define a map $d : {\cal A} \mapsto {\cal A}^\star$ as follows:

Let $a \in A$. Let $c_a : X \mapsto \{a\}$ be the constant map with
value $a$. Set $d(a) = [c_a]$.

As usual, it is an immediate consequence of the {\L}os theorem that
$d$ is an elementary embedding.

\subsubsection{The automorphism $k$}

We define a bijection $s: X \mapsto X$ by $s(x)(n) = x(n+1)$ [for $n
\in \mbox{\bf Z}$].

We define a map $K: {\cal F} \mapsto {\cal F}$ by $K(f)(x) = f(s(x))$.
[The routine check that $K(f)$ is indeed in ${\cal F}$ is left to the
reader.] 

The reader should also verify the more precise statement that if $f$ is
supported by $\{s_1, \ldots, s_n\}$ then $K(f)$ is supported by $\{s_1
+1, \ldots, s_n + 1 \}$.

\begin{lemma}
Let $\varphi(v_1, \ldots, v_n)$ be a formula of the language ${\cal
L}$ appropriate to ${\cal A}$ having the indicated free variables. Let
$f_1, \ldots, f_n \in {\cal F}$.

Then the sentence
\[ \varphi([f_1], \ldots, [f_n]) \iff \varphi([K(f_1)], \ldots,
[K(f_n)])\]
holds in ${\cal A}^\star$.
\end{lemma}

Proof: It is clear that if an element of ${\cal F}$ is supported by a
finite set $s$, then it is supported also by any larger finite set
$s'$. Hence, we may assume that $f_1, \ldots, f_n$ are all supported
by $[a,b]$. [Here, $a$ and $b$ are integers with $a \leq b$. $[a,b] =
\{x \in \mbox{\bf Z}\mid a \leq x \leq b\}$.]

Applying the main result of section~\ref{partition}, we see that there
is an integer $m$ and a function $G: \kappa \mapsto \kappa$ [lying in
$L_{\xi_\infty}$] and  a truth value $\tau \in \{0,1\}$ such that whenever:
\begin{enumerate}
\item $x(a) < x(a+1) \ldots < x(b)$,
\item $x(i) \in B_m$ for $a \leq i \leq b$, and
\item $x(a), \ldots, x(b)$ are $G$ spread apart
\end{enumerate}
then $\varphi(f_1(x),\ldots, f_n(x))$ receives the value $\tau$ in
${\cal A}$.

Let $s = [a,b+1]$. It follows that whenever $x \in A_{s,m,G}$ then
\[\varphi(f_1(x),\ldots,f_n(x)) \iff \varphi(K(f_1)(x), \ldots,
K(f_n)(x))\]
receives the truth value $1$.

The lemma now follows from the {\L}os theorem.

One immediate consequence of the lemma just proved is that we can
define a map $k:A^\star \mapsto A^\star$ by setting $k([f]) = [K(f)]$.
[I.e., the value of $k([f])$ does not depend on the choice of
representative.]

The other immediate consequence is that the map $k$ is an automorphism
of ${\cal A}^\star$.
\subsubsection{Supports}

Let $x \in A^\star$. We say that a finite subset $s \subseteq
\mbox{\bf Z}$ is a support for $x$ if it is a support for some $f \in
{\cal F}$ with $x = [f]$.

We say that $s$ is a block support for $x$ if $s$ is a support for $x$
of the form $[a,b]$ [where $a$ and $b$ are members of {\bf Z} with $a
\leq b$]. Of course every $x \in A^\star$ has a block support.

\begin{lemma}
Let $x \in A^\star$, Then $x$ has a {\em minimum} block support $s_0$
which is contained in every other block support for $x$.
\end{lemma}

Remark: The lemma remains true if the word ``block'' is deleted
throughout. We shall not prove this stronger result.

Proof: Let $s_0$ be a block support for $x$ of minimum cardinality. We
have to show that it is contained in every other block support for
$x$. Suppose not, then there is another block support for $x$ which
neither contains or is contained in $s_0$. 

So there is some element of $s_0$ which is not contained in $s_1$. If
we shring $s_1$ this will continue to be true. Thus we may suppose
that $s_1$ is {\em minimal} in the sense that no proper subset of it is
a block support for $x$. [Obviously $s_0$ is minimal as well.]

Since $s_0$ and $s_1$ are blocks, neither of which is contained in the
other, the least elements of the two blocks must be different.
Interchanging the two blocks, if necessary, we may assume that the
least element of $s_0$ is less than the least element of $s_1$. It
follows easily that also the largest element of $s_0$ is less than the
largest element of $s_1$.

If we performed this interchange, we no longer know that $s_0$ is of
minimum cardinality. But we still have that $s_0$ is minimal [since
before the interchange both $s_0$ and $s_1$ were minimal].

We let $s_0 = [a_0,b_0]$ and $s_1 = [a_1,b_1]$. Let $n_i = card(s_i)$.
Then there are functions $F_i: \kappa^{n_i} \mapsto A$ in
$L_{\xi_\infty}$ such that
defining $f_i: X \mapsto A$ by
\[ f_i(x) = F_i(x(a_i), \ldots, x(b_i))\]
we have $x = [f_0] = [f_1]$.

It follows [by the {\L}os theorem] that there is an integer $m$ and a
function $G:\kappa \mapsto \kappa$ such that whenever 
\begin{enumerate}
\item $x(a_0) < \ldots < x(b_1)$;
\item $x(i) \in B_m$ [for $a_0 \leq i \leq b_1$];
\item $x(a_0), \ldots, x(b_1)$ are $G$ spread apart
\end{enumerate}
then $F_0(x(a_0), \ldots,x(b_0)) = F_1(x(a_1), \ldots, x(b_1))$.

\begin{proposition}
Let $\theta_0, \ldots, \theta_{n_0}$ be an increasing sequence of
ordinals from $B_m$ which are $G$ spread apart. Then $F_0(\theta_0,
\theta_2. \ldots, \theta_{n_0}) = F_0(\theta_1, \ldots, \theta_{n_0})$.
\end{proposition}

Proof: Let $n_2 = card([a_0,b_1])$. We can enlarge the sequence of
thetas to $\theta_0, \ldots, \theta_{n_2}$ so that it is still true
that the sequence of thetas is increasing, consists of elements of
$B_m$, and is $G$ spread apart. 

Let $r = n_2 -n_1 + 1$. Note that our hypotheses on $s_0$ and $s_1$
imply that $r \geq 2$. It follows that 
\[F_0(\theta_0,\theta_2, \ldots, \theta_{n_0}) = F_1(\theta_r, \ldots,
\theta_{n_2}) = F_0(\theta_1, \theta_2, \ldots, \theta_{n_0})\]

The proposition is now clear.

Let now $\theta^\star$ be the least element of $B_m$ which is greater
than $G(0)$. Define $G':\kappa \mapsto \kappa$ so that 
\begin{enumerate}
\item $G'(0) = \max(G(0), \theta^\star + 1, G(\theta^\star))$;
\item $G'(\alpha) = G(\alpha)$ for $\alpha > 0$.
\end{enumerate}

Then if $\theta_1, \ldots, \theta_{n_0}$ is $G'$ spread apart, then
the sequence $\theta^\star, \theta_1, \ldots, \theta_{n_0}$ is
increasing and $G$ spread apart.

Define a function $H: \kappa^{n_0-1} \mapsto A$ by $H(\theta_2,
\ldots, \theta_{n_0}) = F_0(\theta^\star, \theta_2, \ldots,
\theta_{n_0})$. Define $h: X \mapsto A$ by $h(x) = H(x(a_0 +1),
\ldots, x(b_0))$. Then it follows readily from the preceding that
$[f_1] = [h]$. But this contradicts the fact that $s_0$ is a minimal
support for $x$. The lemma is proved.

\subsection{Proof of the main lemma}

We turn now th the proof of Lemma~\ref{main_lemma}. We let ${\cal N}$
and $S$ be as in the statement of that lemma. We choose our set $A_0$
[the input to the ``length $\omega$ construction''] to encode both of
these objects. There is no difficulty doing this since $N$, the
underlying set of ${\cal N}$, has cardinality $\kappa$.

The model to which we will apply our ultrapower construction is
${\cal A} = \langle N\mid \in_N, j\rangle$.

The result of the ultrapower construction is ${\cal A}^\star = \langle
N^\star \mid \in_{N^\star}, j^\star\rangle$.

$N^\star$ will be the underlying set of the model ${\cal N}'$ that we
are constructing. The elementary embedding $\pi$ will just be the
diagonal map $d$ of the preceding subsection.

We have to verify that $\pi$ is not onto. This is easy. Let $i$ be the
given embedding of $L_\kappa$ into $N$. Define $H:X \mapsto N$ by
$H(x) = i(x(0))$. Then it is easily checked that $[H]$ is not in the
range of $\pi$.

The map $i'$ will be the composition $\pi i$. The commutativity of the
diagram of Figure
1 is immediate, and since both $\pi$ and $i$ are elementary embeddings
of models of set theory, so is $i'$.

We have to verify that the range of $i'$ is an initial segment of
$N'$. Suppose not. Then there is an element ${i'}(x)$ and an element $[f]
\in N^\star$ such that $[f] \in_{N^\star} i'(x)$ but $[f]$ is not of
the form $i'(y)$.

Now $f(x)$ = $F(x(a), \ldots, x(b))$ for some $F \in L_{\xi_\infty}$
and $a \leq b$ in {\bf Z}.
By {\L}os, on some set of measure $1$, $f$ takes values in the set
$\{i(y)\mid y \in x\}$ of cardinality less than $\kappa$. It follows
readily from lemma~\ref{partition_lemma} that on some smaller set of
measure $1$, $f$ is constant. This shows that $[f]$ has the form
$i'(y)$.

It follows readily from the fact that $k$ is an automorphism of ${\cal
A}^\star$ and the fact that $j^\star$ is one of the components of
${\cal A}^\star$ that $j^\star$ and  $k$ commute. 

We set $j' = kj^\star$. Since both $k$ and $j$ are automorphisms of
$\langle N^\star \mid \in_{N^\star}\rangle$, so is $j'$.

We must check that the diagram of Figure 2 commutes. Let $a \in N$. We
have to show that $j'(\pi(a)) = \pi(j(a))$. We write $j' = j^\star k$.
Now $\pi$ is just the diagonal map, and we know that elements in the
range of the diagonal map are fixed by $k$. So we are reduced to
proving that $j^\star(\pi(a)) = \pi(j(a))$. But this is clear from the
following facts:
\begin{enumerate}
\item $\pi$ is an elementary embedding from ${\cal A}$ to ${\cal
A}^\star$.
\item $j$ is one of the basic operations of ${\cal A}$.
\item $j^\star$ is the corresponding operation of ${\cal A}^\star$.
\end{enumerate}

Let $x$ be an element of $N^\star$ fixed by $j'$. We must show that
$x$ is in the range of $i'$.

We first show that $x$ is in the range of $\pi$. This amounts to
showing that $x$ has support $\emptyset$. Suppose not toward a
contradiction.

Let $[a,b]$ be the minimum block support for $x$. Let $x = [f]$ with
$f(x) = F(x(a), \ldots, x(b))$ and $F \in L_{\xi_\infty}$. Then $j'(x)$
is represented by the map $\{x \mapsto j(F(x(a+1), \ldots, x(b+1)))$.
Since $j'(x) = x$ and $[a,b]$ is the minimal block support for $x$, we
conclude that $[a,b] \subseteq [a+1,b+1]$ which is absurd. The upshot
is that the element $x$ has empty support and so is in the range of
$\pi$. Say $x = \pi(y)$.

Now \[\pi(y) =x = j'(x) = j'(\pi(y)) = \pi(j(y))\]
Since $\pi$ is injective, we conclude that $y = j(y)$. Since ${\cal
N}$ is an $A$-model, we conclude that $y = i(z)$ for some $z$. But
then $x = \pi(i(z)) = i'(z)$ as desired.

We have now checked that ${\cal N}'$ is an $A$-model, and that
$\pi$ is an $A$-model map which is not onto. The only remaining point
to check is that $S$ is coded in ${\cal N}'$.

Consider the map $h: X \mapsto N$ given by $h(x) = i(S \cap L_{x(0)})$.
It is easy to check that the element $[h]$ codes $S$ in ${\cal N}'$.
The proof of Lemma~\ref{main_lemma} and hence the proof of
Theorem~\ref{main_theorem} is complete.

\section{Postscript}
\label{sec:postscript}

After completing this draft, I discovered a more conceptual way to
think about the ``ultrapower'' construction of
section~\ref{ultrapower}.

In the next draft of this paper, I will incorporate this
improvement. In the meantime, I indicate the main ideas in this
postscript.

\begin{enumerate}
\item One can view the length $\omega$ construction of
  section~\ref{length-omega} as defining a measure, $\nu$, on the
  subsets of $\kappa$ lying in $L_{\xi_\infty}$. Given $B \subseteq
  \kappa$ with $B \in L_{\xi_\infty}$, the construction insures that
  for some $n$, either a tail of $B_n$ is included in $B$ {\em or} a
  tail of $B_n$ is included in $\kappa - B$. In the former case, we
  set $\nu(B)=1$; in the latter case, we set $\nu(B)=0$.

\item By slightly modifiying the construction, we can arrange that for
  any $G: \kappa \mapsto \kappa$ with $G \in L_{\xi_\infty}$, a tail
  of $B_n$ is $G$ spread apart for all sufficiently large $n$.

This is easy to achieve by a suitable thinning of the output of the
``basic module'' of section~\ref{basic-module}. We take $B_{n+1}$ to
be this thinned down set rather than the output of the basic module
applied to $B_n$.

\item The result of this is that the measure $\nu$ will enjoy the
  following partition property. [Cf. the proof in
  section~\ref{partition}.]

If $F: [\kappa]^n \mapsto \gamma$ [with $\gamma < \kappa$ and  $F
\in L_{\xi_\infty}$] then there is a set $B \in L_{\xi_\infty}$ with
$\nu(B) = 1$ such that $F$ is constant on $[B]^n$.

\item One can now go through the construction of
  section~\ref{ultrapower} much as before. But it looks much more
  familiar, resembling the usual iterated ultrapowers with respect to
  a measurable cardinal [with, however, the set {\bf Z} indexing the
  iterations].

\item There is one point to be cautious about. Unlike the usual
  measures constructed from weakly compact cardinals, there is no
  reason to suppose that the ultrapower of $\kappa$ with respect to
  $\nu$ [using only functions $F:\kappa \mapsto \kappa$ lying in
  $L_{\xi_\infty}$] is well-founded or even that there is a least
  non-constant function.
\end{enumerate}
\end{document}